\documentclass[reqno,11pt]{amsart}

\usepackage{hyperref,tikz,amssymb,enumerate}
\usepackage{eucal}

\definecolor{mhcblue}{HTML}{0066CC}

{\theoremstyle{plain}
	\newtheorem{theorem}{Theorem}
	\newtheorem*{theorem*}{Theorem}
	\newtheorem{lemma}[theorem]{Lemma}
	
	\newtheorem{corollary}[theorem]{Corollary}

}
{\theoremstyle{definition}
	\newtheorem{definition}[theorem]{Definition}
	\newtheorem{remark}[theorem]{Remark}
	\newtheorem*{remark*}{Remark}
	
}

\newcommand{\NC}{\mathrm{NC}}
\newcommand{\NCN}{\mathrm{NCN}}
\def\dcr{\mathrm{dcr}}
\def\dne{\mathrm{dne}}
\def\red#1{\mathrm{std}(#1)}

\newcommand{\oeislink}[1]{\href{https://oeis.org/#1}{#1}}

\title[Simple bijection for $k$-noncrossing partitions]{A simple bijection for enhanced, classical, and 2-distant $k$-noncrossing partitions}

\author[Gil]{Juan B. Gil}
\address{Penn State Altoona, 3000 Ivyside Park, Altoona, PA 16601, USA} 
\email{jgil@psu.edu}

\author[Tirrell]{Jordan O. Tirrell}
\address{Department of Mathematics and Computer Science, Washington College, Chestertown, MD 21620, USA}
\email{jtirrell2@washcoll.edu}

\thanks{The authors would like to thank the Department of Mathematics and Statistics and the Hutchcroft fund at Mount Holyoke College for their support.}


\begin{document}
	
\begin{abstract}
In this note, we give a simple extension map from partitions of subsets of $[n]$ to partitions of $[n+1]$, which sends $\delta$-distant $k$-crossings to $(\delta+1)$-distant $k$-crossings (and similarly for nestings). This map provides a combinatorial proof of the fact that the numbers of enhanced, classical, and $2$-distant $k$-noncrossing partitions are each related to the next via the binomial transform. Our work resolves a recent conjecture of Zhicong Lin and generalizes earlier reduction identities for partitions.
\end{abstract}

\maketitle
	
\section{Introduction}

Given a partition $\pi$ of a set of integers and $k\ge 1$, the \textit{arc digram} of $\pi$ is obtained by drawing an arc between each pair of integers that appears consecutively in the same block of $\pi$. For example, the partition $\pi=\{148,26,3,57\}$ can be represented as 
	
	\medskip
	\begin{center}	
		\begin{tikzpicture}[scale=0.45]
		\tikzstyle{node1}=[circle, inner sep=2, fill=black]
		\tikzstyle{arc2}=[draw, line width=1.5,color=black]
		\foreach \x in {0,...,7}{\node[node1] at (2*\x,0) {};};
		\foreach \x in {1,...,8}{\node[below=2pt] at (2*\x-2,0) {\scriptsize \x};};
		\draw[arc2] (0,0) parabola bend (3,1.6) (6,0);		
		\draw[arc2] (2,0) parabola bend (6,1.8) (10,0);		
		\draw[arc2] (6,0) parabola bend (10,1.8) (14,0);		
		\draw[arc2] (8,0) parabola bend (10,1) (12,0);		
		\end{tikzpicture}
	\end{center}
	and we say that $1246$, $2468$ and $2567$ are \textit{crossings},  and $4578$ is a \textit{nesting}. Additionally, we say $148$ is an \textit{enhanced crossing} and $134$ is an \textit{enhanced nesting}. 
	
	Precisely, an \textit{enhanced $k$-crossing} is a sequence 
\begin{equation*}
	a_1<a_2<\cdots<a_k\le b_1<b_2<\cdots<b_k 
\end{equation*} 
such that there is an arc between each pair $a_i,b_i$, and we consider singleton blocks to have trivial arcs. A \textit{classical $k$-crossing} additionally requires that $a_k<b_1$, and in general a \textit{$\delta$-distant $k$-crossing} requires that the \textit{distance} $b_1-a_k$ is at least $\delta$. Thus, the enhanced crossings correspond to $\delta=0$ and the classical crossings correspond to $\delta=1$. Note that an enhanced 1-crossing is simply an arc, and a classical 1-crossing is a nontrivial arc. An \textit{enhanced/classical/$\delta$-distant $k$-nesting} is defined similarly with an arc between each pair $a_i,b_{k+1-i}$. Returning to our example above, there is one enhanced 3-crossing $12468$, and no enhanced 3-nesting. Of the classical crossings and nestings, all but one are 2-distant, and none are 3-distant.

	Much of the combinatorial interest in these structures follows the work of Chen, Deng, Du, Stanley, and Yan~\cite{C+07}. They proved a beautiful bijective symmetry between $k$-crossings and $k$-nestings, which D.~Drake \& J.~S.~Kim~\cite{DK09} generalized to $\delta$-distant $k$-crossings and $k$-nestings (see Eq.~\eqref{eq:NCNsym}). Motivation for $\delta$-distant $k$-crossings also comes from the study of RNA structures~\cite{JQR08}.

	Let $\NC_{\delta,k}(n)$ denote the number of \textit{$\delta$-distant $k$-noncrossing partitions} of $[n]=\{1,\ldots,n\}$, that is, partitions of $[n]$ with no $\delta$-distant $k$-crossing. Examples are given in Table~\ref{tab:NCdk_seqs}. In this paper, we give combinatorial proofs of the following binomial transform identities.
	
\begin{theorem}\label{thm:Lin_conj}
	For integers $n\geq 0$ and $k\geq 1$,
	\begin{align}
	\label{eq:d01}
	\NC_{1,k}(n+1)&=\sum_{i=0}^n\binom{n}{i} \NC_{0,k}(i) \;\text{ and}\\
	\label{eq:d12}
	\NC_{2,k}(n+1)&=\sum_{i=0}^n\binom{n}{i} \NC_{1,k}(i).
	\end{align}
\end{theorem}
	
Equation~\eqref{eq:d01} is well-known for $k=2$, and it was recently proven for $k=3$ by Lin~\cite{Lin18}, who conjectured it held for all $k$. In the final version of \cite{Lin18}, Lin announced two different combinatorial proofs of his conjecture in joint work with D.~Kim~\cite{LK18+}. Equation~\eqref{eq:d12} appears to be new for $k>2$. For $k=2$, it is well-known that $\NC_{1,2}(n)=\frac1{n+1}\binom{2n}n$. The fact that $\NC_{2,2}(n)$ is the binomial transform of $\NC_{1,2}(n)$ was shown by Drake \& J.~S.~Kim~\cite{DK09} using Charlier diagrams and generating functions. Several bijections with other objects can be found in \cite{GK10,Kim11}. 

\begin{table}[ht]
\small
\def\R{\rule[-1ex]{0ex}{3.6ex}}
\def\g{\color{gray!50}}
\begin{tabular}{c|cccccccccccc|c}
	$n$& 0 & 1 & 2 & 3 & 4 & 5 & 6 & 7 & 8 & 9 & 10 & 11 & OEIS \cite{OEIS}\\[2pt]
	\hline
	\R $\NC_{0,1}(n)$ &\g 1 & 0 & 0 & 0 & 0 & 0 & 0 & 0 & 0 & 0 & 0 & 0 &  \\
	$\NC_{0,2}(n)$ &\g 1 &\g 1 &\g 2 & 4 & 9 & 21 & 51 & 127 & 323 & 835 & 2188 & 5798 & \oeislink{A001006} \\
	$\NC_{0,3}(n)$ &\g 1 &\g 1 &\g 2 &\g 5 &\g 15 & 51 & 191 & 772 & 3320 & 15032 & 71084 & 348889 & \oeislink{A108307} \\
	$\NC_{0,4}(n)$ &\g 1 &\g 1 &\g 2 &\g 5 &\g 15 &\g 52 &\g 203 & 876 & 4120 & 20883 & 113034 & 648410 & \oeislink{A192855} \\
	$\NC_{0,5}(n)$ &\g 1 &\g 1 &\g 2 &\g 5 &\g 15 &\g 52 &\g 203 &\g 877 &\g 4140 & 21146 & 115945 & 678012 & \oeislink{A192865} \\[2pt]
	\hline
	\R $\NC_{1,1}(n)$ &\g 1 &\g 1 & 1 & 1 & 1 & 1 & 1 & 1 & 1 & 1 & 1 & 1 & \\
	$\NC_{1,2}(n)$ &\g 1 &\g 1 &\g 2 &\g 5 & 14 & 42 & 132 & 429 & 1430 & 4862 & 16796 & 58786 & \oeislink{A000108} \\
	$\NC_{1,3}(n)$ &\g 1 &\g 1 &\g 2 &\g 5 &\g 15 &\g 52 & 202 & 859 & 3930 & 19095 & 97566 & 520257 & \oeislink{A108304} \\
	$\NC_{1,4}(n)$ &\g 1 &\g 1 &\g 2 &\g 5 &\g 15 &\g 52 &\g 203 &\g 877 & 4139 & 21119 & 115495 & 671969 & \oeislink{A108305} \\
	$\NC_{1,5}(n)$ &\g 1 &\g 1 &\g 2 &\g 5 &\g 15 &\g 52 &\g 203 &\g 877 &\g 4140 &\g 21147 & 115974 & 678530 & \oeislink{A192126} \\[2pt]
	\hline
	\R $\NC_{2,1}(n)$ &\g 1 &\g 1 &\g 2 & 4 & 8 & 16 & 32 & 64 & 128 & 256 & 512 & 1024 & \oeislink{A000079} \\
	$\NC_{2,2}(n)$ &\g 1 &\g 1 &\g 2 &\g 5 &\g 15 & 51 & 188 & 731 & 2950 & 12235 & 51822 & 223191 & \oeislink{A007317} \\
	$\NC_{2,3}(n)$ &\g 1 &\g 1 &\g 2 &\g 5 &\g 15 &\g 52 &\g 203 & 876 & 4115 & 20765 & 111301 & 627821 & \oeislink{A366774} \\
	$\NC_{2,4}(n)$ &\g 1 &\g 1 &\g 2 &\g 5 &\g 15 &\g 52 &\g 203 &\g 877 &\g 4140 & 21146 & 115938 & 677765 & \oeislink{A366775} \\
	$\NC_{2,5}(n)$ &\g 1 &\g 1 &\g 2 &\g 5 &\g 15 &\g 52 &\g 203 &\g 877 &\g 4140 &\g 21147 &\g 115975 & 678569 & \oeislink{A366776} \\[2pt]
	\hline
	\R $B(n)$ & 1 & 1 & 2 & 5 & 15 & 52 & 203 & 877 & 4140 & 21147 & 115975 & 678570 & \oeislink{A000110}
\end{tabular}
\bigskip
\caption{Sequences $\NC_{\delta,k}(n)$ for small $k$ and $\delta=0,1,2$.}
\label{tab:NCdk_seqs}
\end{table}
\vspace{-1ex}
Bousquet-M\'elou \& Xin~\cite{BX06} showed that the sequences  $\NC_{0,3}(n)$ and $\NC_{1,3}(n)$ are P-recursive and gave explicit recurrences. Mishna \& Yen~\cite{MY13} used generating trees to find functional equations for $\NC_{1,k}(n)$ when $k>3$. Burrill, Elizalde, Mishna, and Yen~\cite{BEMY16} did the same for $\NC_{0,k}(n)$. Theorem~\ref{thm:Lin_conj} provides a simple way to connect the enumerations of $\NC_{0,k}(n)$, $\NC_{1,k}(n)$, and $\NC_{2,k}(n)$. In particular, as Lin~\cite{Lin18} observed, the binomial transform preserves D-finiteness (see~\cite[Theorem 6.4.10]{Sta99}).
\begin{corollary}
 The D-finiteness of the ordinary generating functions for $\NC_{0,k}(n)$, $\NC_{1,k}(n)$, and $\NC_{2,k}(n)$ are the same.\footnote{For $k>3$, it is conjectured~\cite{BX06} that the generating function of $\NC_{1,k}(n)$ is not D-finite.}
\end{corollary}

Unfortunately, the sequence $\NC_{3,k}(n+1)$ is not the binomial transform of the sequence $\NC_{2,k}(n)$, and in general, not much is known about $\NC_{\delta,k}(n)$ for $\delta\ge 3$ (the generating function for $\NC_{3,2}(n)$ is given in \cite{Kim11-fr}).

Let $\NCN_{\delta,k,\varepsilon,j}(n)$ denote the number of $\delta$-distant $k$-noncrossing partitions of $[n]$ which also contain no $\varepsilon$-distant $j$-nesting. By the celebrated symmetry of Chen et al.~\cite{C+07}, generalized by Drake \& Kim~\cite{DK09}, we have
\begin{equation}\label{eq:NCNsym}
	\NCN_{\delta,k,\varepsilon,j}(n)=\NCN_{\varepsilon,j,\delta,k}(n).
\end{equation}

We will give a combinatorial proof of the following refinement of Theorem~\ref{thm:Lin_conj}.
\begin{theorem}\label{thm:NCN} For integers $n\ge 0$, $j,k\ge 1$, and $\varepsilon,\delta\in\{0,1\}$,
	\begin{equation}
	\NCN_{\delta+1,k,\varepsilon+1,j}(n+1)=\sum_{i=0}^n\binom ni\NCN_{\delta,k,\varepsilon,j}(i)
	\end{equation}
\end{theorem}

The generating function for $\NCN_{1,k,1,j}(n)$ is not simply D-finite but in fact rational for all $k,j\ge 1$. This was shown by Marberg~\cite{Mar13} using a bijection with walks on certain multigraphs (and he actually proved this in the more general setting of colored partitions). By Theorem~\ref{thm:NCN}, the generating functions for $\NCN_{0,k,0,j}(n)$ and $\NCN_{2,k,2,j}(n)$ are also rational.

\section{Main bijection}

	It is known that the sequence $B(n)$ of Bell numbers, which counts the number of partitions of $[n]$, is an eigensequence of the binomial transform. That is, $B(n+1) = \sum_{i=0}^n\binom{n}{i} B(i)$. The usual combinatorial proof of this fact uses a bijection from partitions of subsets of $[n]$ to partitions of $[n+1]$, defined by adding a new block containing $n+1$ together with every element of $[n]$ not already contained in a block. For example, if $n=9$, this bijection takes $\{1479,25,6\}$ to $\{1479,25,6,38\text{T}\}$, where T stands for 10. 
We now give a different bijection which behaves well with respect to crossings and distance.
	
\begin{definition}
	For any $n\ge 0$, and for any partition $\pi$ of some subset of $[n]$, we define the partition $\hat{\pi}$ of $[n+1]$ as follows. For each pair of vertices $v<w$ appearing consecutively in the same block of $\pi$, place $v$ and $w+1$ in the same block of $\hat{\pi}$. For each singleton $u$ of $\pi$, place $u$ and $u+1$ in the same block of $\hat{\pi}$. At last, place any remaining vertices in $[n+1]$ as singletons of $\hat{\pi}$. We will refer to $\pi\mapsto\hat{\pi}$ as the \textit{extension map}.
\end{definition}

	For example, if $n=9$ and $\pi=\{1479,25,6\}$, then $\hat\pi=\{15,267\text{T},3,48,9\}$. An elegant geometric description of this map can be obtained by extending the lines in the arc diagram of $\pi$ to arrive at the arc diagram of $\hat{\pi}$ (see Figure~\ref{fig:vis_bij}), where singletons are considered to have trivial arcs.

An earlier description of our map $\pi\mapsto\hat{\pi}$ followed a suggestion by Lin~\cite{Lin18} and used growth diagrams \`a la Krattenthaler~\cite{Kratt06} (as in work by Kasraoui~\cite{Kas09} and by Yan~\cite{Yan18}).

\medskip
\begin{figure}[ht]
	\usetikzlibrary{shapes}
	\begin{tikzpicture}[scale=0.65]
	\tikzstyle{node0}=[circle, inner sep=2, fill=mhcblue]
	\tikzstyle{node1}=[circle, inner sep=2, fill=black]
	\tikzstyle{arc0}=[draw, line width=1.4, dotted, color=mhcblue]
	\tikzstyle{arc1}=[draw, line width=1.5, solid, color=black]
	\draw[step=1,lightgray,thin] (0,0) grid (18,3);
	\foreach \x in {0,...,9}{\node[node0] at (2*\x,0) {};};
	\foreach \x in {1,...,10}{\node[color=mhcblue,below=2pt] at (2*\x-2,0) {\scriptsize \x};};
	\draw[arc0] (0.05,0) parabola bend (4,2.4) (7.95,0);		
	\draw[arc0] (2.05,0) parabola bend (6,2.4) (9.95,0);		
	\draw[arc0] (6.05,0) parabola bend (10,2.4) (13.95,0);		
	\draw[arc0] (10.05,0) parabola bend (11,1) (11.95,0);		
	\draw[arc0] (12.05,0) parabola bend (15,1.9) (17.95,0);		
	\foreach \x in {0,1,3,4,5,6,8}{\node[node1] at (1+2*\x,1) {};};	
	\foreach \x in {1,2,4,5,6,7,9}{\node[below=1pt] at (2*\x-1,1) {\scriptsize \x};};	
	\draw[arc1] (0.98,1) parabola bend (4,2.4) (7.02,1);		
	\draw[arc1] (6.98,1) parabola bend (10,2.4) (13.02,1);		
	\draw[arc1] (12.98,1) parabola bend (15,1.9) (17.02,1);		
	\draw[arc1] (2.98,1) parabola bend (6,2.4) (9.02,1);		
	\end{tikzpicture}
	\caption{The arc diagram for $\pi=\{1479,25,6\}$ extends to $\hat\pi=\{15,267\text{T},3,48,9\}$.}
	\label{fig:vis_bij}
\end{figure}
	
	It is straightforward to confirm that $\pi\mapsto\hat{\pi}$ is a bijection from partitions of subsets of $[n]$ to partitions of $[n+1]$, and that it increases the distance of $k$-crossings (and $k$-nestings) by one. 
In particular, $(\delta+1)$-distant $k$-crossings in $\hat{\pi}$ correspond to $\delta$-distant $k$-crossings in $\pi$, which for $\delta\in\{0,1\}$ correspond to $\delta$-distant $k$-crossings in the \textit{standardization} $\red\pi$.%
\footnote{If $\pi$ is a partition of $X\subseteq[n]$, then $\red{\pi}$ is obtained by relabeling the elements of $\pi$ to $1,2,\ldots,|X|$, while preserving their order.}
However, for $\delta\ge 2$, the distance of a $k$-crossing in $\pi$ may not be the same as in $\red{\pi}$. For example, the $3$-distant crossing $1<2<5<6$ in $\hat\pi=\{15,26,3,4\}$ is the image of the $2$-distant crossing $1<2<4<5$ in $\pi=\{14,25\}$, which corresponds to the $1$-distant crossing $1<2<3<4$ in $\red\pi =\{13,24\}$. This explains why $\NC_{\delta+1,k}(n+1)$ is not the binomial transform of the sequence $\NC_{\delta,k}(n)$ for $\delta\ge 2$.

\medskip
The following lemma follows immediately from our definitions.
\begin{lemma}\label{lem:main}
The extension map $\pi\mapsto\hat{\pi}$ satisfies the following properties, where $S$ is any sequence $(a_1,\ldots,a_k,b_1,\ldots,b_k)$, $\widehat S=(a_1,\ldots,a_k,b_1+1,\ldots,b_k+1)$, $\delta\ge0$, and $k\ge 1$.
\begin{enumerate}[$(i)$]
	\item The sequence $S$ is a $k$-crossing of distance $\delta$ in $\pi$ if and only if $\widehat S$ is a $k$-crossing of distance $\delta+1$ in $\hat\pi$.
	\item The sequence $S$ is a $\delta$-distant $k$-crossing after the relabeling from $\pi$ to $\red\pi$ only if $\widehat S$ is a $(\delta+1)$-distant $k$-crossing in $\hat{\pi}$.
	\item The sequence $S$ is an enhanced (resp.~classical) $k$-crossing after the relabeling from $\pi$ to $\red\pi$ if and only if $\widehat S$ is a classical (resp.~2-distant) $k$-crossing in $\hat{\pi}$.
\end{enumerate}
Moreover, analogous properties hold for nestings.
\end{lemma}

In particular, if $\dcr_{\delta,k}(\pi)$ (resp.~$\dne_{\delta,k}(\pi)$) denotes the number of $\delta$-distant $k$-crossings (resp.~$k$-nestings) of a given partition $\pi$, then
\begin{align}
	\dcr_{\delta,k}(\red\pi)&\le \dcr_{\delta,k}(\pi) = \dcr_{\delta+1,k}(\hat{\pi}),\;\text{ and}\\
	\dne_{\delta,k}(\red\pi)&\le \dne_{\delta,k}(\pi) = \dne_{\delta+1,k}(\hat{\pi})\;\text{ for all }\delta\ge0,\, k\ge1,
\end{align}
with equality when $\delta\in\{0,1\}$. See Table~\ref{tab:dcr_dne} for examples.

\begin{table}[ht!]
\def\R{\rule[-1ex]{0ex}{3.6ex}}
	\begin{tabular}{c|ccc|ccc|ccc}
	&$\dcr_{0,1}$&$\dcr_{1,1}$&$\dcr_{2,1}$&$\dcr_{0,2}$&$\dcr_{1,2}$&$\dcr_{2,2}$&$\dcr_{0,3}$&$\dcr_{1,3}$&$\dcr_{2,3}$\\[1pt]\hline
	\R $\red\pi$&5&4&3&4&2&0&1&0&0\\
	$\pi$&5&4&4&4&2&1&1&0&0\\
	$\hat\pi$&7&5&4&6&4&2&1&1&0%
	\vspace{1em}\\
	&$\dne_{0,1}$&$\dne_{1,1}$&$\dne_{2,1}$&$\dne_{0,2}$&$\dne_{1,2}$&$\dne_{2,2}$&$\dne_{0,3}$&$\dne_{1,3}$&$\dne_{2,3}$\\[1pt]\hline
	\R $\red\pi$&5&4&3&1&0&0&0&0&0\\
	$\pi$&5&4&4&1&0&0&0&0&0\\
	$\hat\pi$&7&5&4&4&1&0&0&0&0\\[2pt]
	\end{tabular}
\bigskip
\caption{Values of $\dcr_{\delta,k}$ and $\dne_{\delta,k}$ for $\pi=\{1479,25,6\}$.}
\label{tab:dcr_dne}
\end{table}

Lemma~\ref{lem:main} immediately gives us the following generalization of Theorems~\ref{thm:Lin_conj} and~\ref{thm:NCN}.

\begin{theorem} 
For any $n\ge 0$ and given sequences $\alpha, \alpha', \beta$, and $\beta'$, let $P_{\delta}(n,\alpha,\alpha',\beta,\beta')$ denote the number of partitions of $[n]$ with $\dcr_{\delta,k}=\alpha_k$, $\dcr_{\delta+1,k}=\alpha'_k$, $\dne_{\delta,k}=\beta_k$, $\dne_{\delta+1,k}=\beta'_k$ for all $k\ge 1$. Then
\begin{equation}
	P_1(n+1,\alpha,\alpha',\beta,\beta')=\sum_{i=0}^n \binom{n}{i}P_0(i,\alpha,\alpha',\beta,\beta').
\end{equation}
\end{theorem}

Summing over all $\alpha,\alpha',\beta,\beta'$ with either $\alpha_k=0$ or $\alpha'_k=0$ gives Equation \eqref{eq:d01} or \eqref{eq:d12} of Theorem~\ref{thm:Lin_conj}. 
Summing over all $\alpha,\alpha',\beta,\beta'$ with either $\alpha_k=0$ or $\alpha'_k=0$, and either $\beta_j=0$ or $\beta'_j=0$ gives the different $\delta,\varepsilon\in\{0,1\}$ cases of Theorem~\ref{thm:NCN}.
Summing over all $\alpha,\alpha',\beta,\beta'$ with some fixed $\alpha_1$, and $\alpha_2=0$, we obtain a result equivalent to~\cite[Theorem 3.2]{CDD05}. In fact, the ``reduction algorithm" used in \cite{CDD05} is essentially $\hat{\pi}\mapsto\pi$, the inverse of our extension map. 

\begin{remark}\label{rem:dcrprops}
Note that placing restrictions on $\dcr_{\delta,k}$ and $\dne_{\delta,k}$ for $\delta\ge 0$ and $k\ge1$ allows us to control many important properties of partitions. Below we highlight a few key features. Similar statements hold for nestings.
\begin{enumerate}[\quad$\triangleright$]
	\item $\dcr_{\delta,k}(\pi)-\dcr_{\delta+1,k}(\pi)$ is the number of $k$-crossings of distance exactly $\delta$.
	\item $\dcr_{0,1}(\pi)-\dcr_{1,1}(\pi)=\dne_{0,1}(\pi)-\dne_{1,1}(\pi)$ is simply the number of singleton blocks, whose elements are both minimal and maximal in their block.
	\item $\dcr_{0,2}(\pi)-\dcr_{1,2}(\pi)$ is the number of \textit{transients}, elements which are neither minimal nor maximal in their block.
	\item $\dcr_{1,1}(\pi)=\dne_{1,1}(\pi)$ is the number of edges in the arc diagram of $\pi$, which equals $n-k$ when $\pi$ is a partition of $[n]$ with $k$ blocks.
	\item $\dcr_{\delta,k}(\pi)=0$ if and only if $\pi$ is $\delta$-distant $k$-noncrossing.
	\item $\dcr_{m,1}(\pi)-\dcr_{1,1}(\pi)=\dne_{m,1}(\pi)-\dne_{1,1}(\pi)$ is the number of consecutive pairs in blocks at distance less than $m$. When this is zero, $\pi$ is said to be $m$-regular.
\end{enumerate}
\end{remark}

\section{Further remarks}

\subsection*{Partitions of $[n]$ to partitions of $[n+1]$}

If we restrict our map $\pi\mapsto\hat\pi$ to partitions of $[n]$ such that $\red{\pi}=\pi$ (i.e.~there is no $i\in[n]$ which is not contained in some block of $\pi$), we obtain a bijection between partitions $\pi$ of $[n]$ and partitions $\hat{\pi}$ of $[n+1]$ having no $i\in[n]$ such that $i$ is maximal within its block in $\hat\pi$ and $i+1$ is minimal within its block in $\hat\pi$.

Alternatively, we can map a partition $\pi$ of $[n]$ into a partition of $[n+1]$ by first removing its singletons (denote the modified singleton free partition by $\pi'$), and then applying the extension map to $\pi'$ to obtain a partition $\hat{\pi}'$ of $[n+1]$. For example, if $\pi=\{1478,25,3,6\}$, then $\pi'=\{1478,25\}$ and $\hat{\pi}'=\{15,26,3,48,79\}$. Note that singletons only occur in enhanced $k$-nestings where the distance is zero, not in classical $k$-nestings. Moreover, singletons occur in enhanced $k$-crossings only when $k=1$ and the distance is zero. Therefore, 
\begin{equation}
  \dne_{\delta,k}(\pi)= \dne_{\delta,k}(\pi') = \dne_{\delta+1,k}(\hat{\pi}') \;\text{ for all } \delta,k\ge1.
\end{equation}
And for all $\delta\ge 0$, $k\ge1$, except the case $\delta=0$ and $k=1$,
\begin{equation}
  \dcr_{\delta,k}(\pi)= \dcr_{\delta,k}(\pi') = \dcr_{\delta+1,k}(\hat{\pi}').
\end{equation}

This gives the following.
\begin{theorem} 
For $n\ge 0$, and for $(\delta\ge 0,\, k>1)$ or $(\delta>0,\, k\ge 1)$, let $\Pi_{\delta,k}(n,\alpha,\beta)$ be the set of partitions of $[n]$ with $\dcr_{\delta,k}=\alpha$ and $\dne_{\delta,k}=\beta$. The map $\pi\mapsto\hat{\pi}'$ gives a bijection between $\Pi_{\delta,k}(n,\alpha,\beta)$ and the set of partitions of $[n+1]$ with $\dcr_{1,1}=\dcr_{1,2}$ (i.e.~2-regular), $\dcr_{\delta+1,k}=\alpha$, and $\dne_{\delta+1,k}=\beta$.
\end{theorem}

\subsection*{Partitions without singleton blocks}
Let $\overline{\NC}_{\delta,k}(n)$ denote the number of $\delta$-distant $k$-noncrossing partitions of $[n]$ with no singleton blocks. Then 
\[ \NC_{\delta,k}(n)=\sum_{i=0}^n\binom ni\overline{\NC}_{\delta,k}(i), \]
and it can be easily checked that
\[ \NC_{\delta,k}(n+1)=\sum_{i=0}^n\binom ni\Big(\overline{\NC}_{\delta,k}(i)+\overline{\NC}_{\delta,k}(i+1)\Big). \]
Using the inverse binomial transform, we then get that Theorem~\ref{thm:Lin_conj} is equivalent to 
\begin{align*}
	\NC_{0,k}(n) &= \overline{\NC}_{1,k}(n) + \overline{\NC}_{1,k}(n+1), \;\text{ and} \\
	\NC_{1,k}(n) &= \overline{\NC}_{2,k}(n) + \overline{\NC}_{2,k}(n+1).
\end{align*}

From these equations, one obtains that $\overline{\NC}_{1,2}(n)$ gives the sequence of Riordan numbers \oeislink{A005043}, and that $\overline{\NC}_{2,2}(n)$ gives the sequence \oeislink{A033297}. The sequence $\overline{\NC}_{0,3}(n)$ has been shown by the second author, Westbury, and Zhang~\cite{TW18+} to be the sequence \oeislink{A059710}, which appears in the representation theory of $G_2$ (see \cite{Wes07}).
	

\end{document}